\documentclass[11pt]{amsart}
\usepackage{amssymb,amsmath}
\usepackage[dvipdfmx]{hyperref}
\usepackage[all]{xy}
\usepackage{color}
\usepackage[T1]{fontenc}

\newtheorem{theorem}{Theorem}[section]
\newtheorem{proposition}[theorem]{Proposition}
\newtheorem{lemma}[theorem]{Lemma}
\newtheorem{corollary}[theorem]{Corollary}

\theoremstyle{definition}
\newtheorem{definition}[theorem]{Definition}

\theoremstyle{remark}

\numberwithin{equation}{section}

\def\tA2{\tilde A_2}
\def\SL{\mathrm{SL}}
\def\St{\mathrm{St}}

\title{Computing Kazhdan constants by semidefinite programming}

\begin{document}

\author{Koji Fujiwara and Yuichi Kabaya}

\maketitle

\begin{abstract}
Kazhdan constants of discrete groups are hard to 
compute and the actual constants are known only for several 
classes of groups. 

By solving a semidefinite programming problem by a computer,  we obtain a lower bound of 
the Kazhdan constant of a discrete group. 
Positive lower bounds imply that the group has property (T).

We study lattices on $\tilde A_2$-buildings in detail.
For $\tilde A_2$-groups, our numerical bounds look identical 
to the known actual constants. That suggests that our approach is 
effective. 

For a family of groups, $G_1, \cdots, G_4$, that are studied
by Ronan, Tits and others, we conjecture the spectral gap
of the Laplacian is $(\sqrt 2-1)^2$ based on our experimental results. 

For $\mathrm{SL}(3,\Bbb Z)$ and $\mathrm{SL}(4,\Bbb Z)$ we obtain lower bounds of the Kazhdan constants,
0.2155 and 0.3285, respectively, which
are better than any other known bounds. 
We also obtain 0.1710 as a lower bound 
of the Kazhdan constant of the Steinberg group
$\St_3(\Bbb Z)$. 

\end{abstract}

\section{Introduction}
Property (T) of a group was introduced by Kazhdan 
and he proved that, for example, $\mathrm{SL}(n,\Bbb Z)$ have 
property (T) if $n \ge 3$.
We say a group is a Kazhdan group if the group has property (T). 
Property (T) plays an important role in many areas of mathematics.

For a Kazhdan group, given a finite generating set, 
a positive constant $\kappa$, called {\it the Kazhdan 
constant} is defined. Lower and upper bounds of the 
Kazhdan constants are known for some examples
of groups, but to compute the value of the Kazhdan
constant is usually hard.

Recently Ozawa \cite{ozawa} found a new characterization of a Kazhdan group among finitely generated groups. A significance of his theorem
is that it gives a semidefinite algorithm to decide if a finitely 
generated group is a Kazhdan group or not.
To be precise, he gives a description of the spectral gap
of the Laplacian, which is closely related to the Kazhdan constant.

Netzer-Thom implemented it as a computer
program using Semidefinite programming. 
They obtain a lower bound, 0.1783,  for the Kazhdan constant 
of  $\mathrm{SL}(3,\Bbb Z)$ with respect to the set of elementary matrices, which 
is much better than any other known bounds. 

Building up on those ideas, we compute lower
bounds of the spectral gaps of the Laplacian 
and the Kazhdan constants in various examples.
As far as we know, among infinite groups, the Kazhdan constants are 
known only for the family of ``$\tilde A_2$-groups''
(Theorem \ref{A2group.bound}).
Those are groups  that act on $\tilde A_2$-buildings properly and cocompactly
in a very special way. 
For those groups, we obtain lower bounds of the Kazhdan constants, 
which are almost identical 
to the actual values (Tables \ref{table:CMSZ2_q=2}
and \ref{table:CMSZ2_q=3}).

We regard this computation as an evidence that our method is 
effective to compute a good lower bound. 

We also point out that this class of groups
is interesting in various ways, for example, 
the equality holds in two basic inequalities
regarding the spectral gaps and the Kazhdan constants
 (see Corollary \ref{gap.achieve}).

We then discuss another class of groups, $G_1, \cdots, G_4$,
which also act on some $\tilde A_2$-buildings. 
This family is found by Ronan \cite{Ronan}, then 
studied from various viewpoints, including by Tits \cite{Tits}
as ``triangles of groups'' (see Section \ref{sec:ronan}).
In particular they are Kazhdan groups, but the Kazhdan constants
are unknown. 
We obtain 0.239146... as a common (numerical) lower
bound of the Kazhdan constants (see Table \ref{table:ronan}), and also 0.171573
as a common (numerical) lower bound of the spectral gaps.
We predict that they are $(\sqrt 2 -1)/\sqrt 3 = 0.239146...$
and $(\sqrt 2-1)^2=0.1715728...$.

Also, we obtain 0.2155 as a lower bound of the 
Kazhdan constant of  $\mathrm{SL}(3,\Bbb Z)$,
which is slightly better(bigger) than the one by 
Netzer-Thom, and also 0.3285  for $\mathrm{SL}(4,\Bbb Z)$,
which is much better than any other known bounds
(see Table \ref{table:SL3Z_SL4Z}).

We were not able to obtain a positive bound for $\mathrm{SL}(5,\Bbb Z)$
because of the lack of the power of the computer. 
Also, unfortunately, we did not find any new examples
of Kazhdan groups by our method. 

The Steinberg group $\St_n(\Bbb Z)$ is very closely 
related to $\SL(n,\Bbb Z)$, which is the quotient of $\St_n(\Bbb Z)$
by $\Bbb Z/2\Bbb Z$.
We obtain 0.1710 as a lower bound 
of the Kazhdan constant of $\St_3(\Bbb Z)$
(see the inequality \ref{steinberg_kappa}).
As far as we know, no concrete lower 
bounds of the Kazhdan constants of the Steinberg groups were known.
Also, our bound is obtained from a finite presentation of 
$\St_3(\Bbb Z)$. In that sense, this is a new proof that $\St_3(\Bbb Z)$
is a Kazhdan group.

Then we turn our attention to finite groups.
All of them are Kazhdan groups, but again, there are a small number of 
cases that the Kazhdan constants are known.
One such example is the family of finite Coxeter groups.
We obtain bounds and they are very close to the 
known values (see Section \ref{sec:coxeter}).
We also obtain lower bounds for complex Coxeter groups,
for which the actual values are not known
(see Section \ref{sec:cpx_ref}).
We believe that our bounds are very close to the actual values. 

At the end, we ask questions that arise naturally from 
our experimental results.

{\bf Acknowledgements}.
We would like to thank Uri Bader, Pierre-Emmanuel Caprace,
Mike Davis, Ian Leary, Pierre Pansu, Narutaka Ozawa
 and Alain Valette. 
We benefit from  Kawakami's paper \cite{Ka}.

 We are supported by Grant-in-Aid for Scientific Research (No. 15H05739).

\section{Kazhdan constants and spectral gaps}
Let $\Gamma$ be a finitely generated group with a finite generating set $S$.
For a Hilbert space $\mathcal{H}$ and 
a unitary representation $\pi : \Gamma \to \mathcal{H}$, set
\[
\kappa( \Gamma, S, \pi) = \inf_{\substack{\xi \in \mathcal{H}, \\ || \xi || = 1}} 
\max_{s \in S} || \pi(s) \xi - \xi ||.
\]
We define the \emph{(optimal) Kazhdan constant} with respect to $S$ by
\[ 
\kappa(\Gamma, S) = \inf \{ \kappa( \Gamma, S, \pi) \mid 
\textrm{$\pi$ has no non-zero invariant vector} \}.
\]
If $\kappa(\Gamma, S) > 0$ for some (then any) finite generating set $S$,
then  the group $\Gamma$ is called a \emph{Kazhdan group}.

From now on, we assume that a generating set $S$ is symmetric.
Let $\mathbb{R}[\Gamma]$ be the group ring of $\Gamma$ with coefficients in $\mathbb{R}$. 
Define the \emph{unnormalized Laplacian} $\Delta$ by
\[
\label{eq:def_of_unnormlized_laplacian}
\Delta = |S| - \sum_{s \in S} s = 
\frac{1}{2} \sum_{s \in S} (1-s)(1-s^{-1}) \in \mathbb{R}[\Gamma].
\]
Define the $*$-operation on $\mathbb{R}[\Gamma]$ by 
$(\sum_{g \in \Gamma} r_g g)^* = \sum_{g \in \Gamma} r_g g^{-1}$ $(r_g \in \mathbb{R})$.
By definition, $\Delta^*=\Delta$.
For any unitary representation $\pi : \Gamma \to U(\mathcal{H})$, 
$\pi( \Delta )$ is an operator on $\mathcal{H}$.
We remark that the normalized Laplacian is defined
as $\Delta/|S|$, which is used in \cite{ozawa}.

\begin{theorem}[Proposition 5.4.5 and Remark 5.4.7 of \cite{BdlHV}]
\label{thm:gap}
Let $\Gamma$ be a discrete group with finite symmetric generating set $S$.
Suppose there exists $\varepsilon > 0$ such that 
for any unitary representation $\pi : \Gamma \to \mathcal{H}$ 
with no non-zero invariant vector, 
\[
\langle \pi(\Delta) \xi, \xi  \rangle \geq \varepsilon \langle \xi, \xi \rangle
\]
for any $\xi \in \mathcal{H}$. Then
\begin{equation}
\label{eq:ek}
\sqrt{\frac{2 \varepsilon}{|S|}} \leq \kappa(\Gamma, S).
\end{equation}
In particular, $\Gamma$ is a Kazhdan group.
\end{theorem}

In some cases, the equality holds in 
the inequality \ref{eq:ek} (see Corollary \ref{gap.achieve}).

We remark \cite{BdlHV} uses the normalized Laplacian, 
therefore the factor $1/|S|$ does not appear in their setting. 
We call $\varepsilon$  {\it a spectral gap}
and the supremum of $\varepsilon$  {\it the spectral gap} of $\Delta$.

\subsection{Ozawa's criterion}
We recall a theorem by Ozawa. 

\begin{theorem}[Ozawa {\cite[Main Theorem]{ozawa}}]
\label{thm:ozawa}
There exist $b_1, \dots, b_n \in \mathbb{R}[\Gamma]$ and $\varepsilon > 0$ such that 
\begin{equation}
\label{eq:ozawa}
\Delta^2 - \varepsilon \Delta = \sum_{i=1}^n b_i^* b_i,
\end{equation}
if and only if $\Gamma$ is a Kazhdan group.

Moreover, if $\Gamma$ is a Kazhdan group then for any rational 
number $\varepsilon' \le \varepsilon$, there exist
positive rational numbers $r_i$ and $b_1, \dots, b_n \in \mathbb{Q}[\Gamma]$
such that 
$\Delta^2 - \varepsilon' \Delta = \sum_{i=1}^n r_i b_i^* b_i$
\end{theorem}

By the spectral mapping theorem, the condition in Theorem \ref{eq:ek} and the condition in Theorem \ref{thm:ozawa}
are equivalent for a same constant $\varepsilon$ (\cite{ozawa}).
\begin{definition}[Spectral gap]
We call $\varepsilon$ and  the supremum of $\varepsilon$, 
denoted by $\varepsilon(\Gamma,S)$, in 
Theorem \ref{thm:ozawa} {\it a spectral gap} and {\it the spectral gap} (in the sense of Ozawa),
respectively. We sometimes write $\varepsilon(\Gamma,S)$
as $\varepsilon$.
\end{definition}

Notice that from the moreover part of the theorem, 
since $\Gamma$ is countable, if $\Gamma$ is a Kazhdan group,
we will be able to find a rational solution with $\varepsilon >0$, which implies
that $\Gamma$ is a Kazhdan group. This is the case
even if a group does not have a solvable word problem since 
we only claim the algorithm is semidefinite (ie, stops
only when there is a positive solution).

  Also a number $\varepsilon$ we obtain gives a lower bound of the spectral gap, which will give a lower bound 
of the Kazhdan constant by the inequality (\ref{eq:ek}).

\subsection{Solving the equation}
In view of this, we  try to find a (rational) solution, $\varepsilon, b_i$, for the equation
(\ref{eq:ozawa})
using a computer to obtain a lower bound of the Kazhdan constant.
Netzer-Thom \cite{netzer-thom} already carried out 
this strategy and found a solution with positive $\varepsilon$
for $\mathrm{SL}(3,\Bbb Z)$. The lower bound they obtain 
for the Kazhdan constant is much better than any known bounds
(see Section \ref{sec:SL3Z})

There is an issue in this strategy. A computer can find only a numerical
solution, so that even if we find a solution  with some $\varepsilon >0$,  maybe the equation does not have any ``exact''
solutions with $\varepsilon >0$.
Netzer-Thom addressed this issue and found 
 a way to certify the existence of an exact 
solution, with an estimate of $\varepsilon$, once we have a ``good'' numerical solution
as the following lemmas show.

 \begin{lemma}[Netzer-Thom \cite{netzer-thom}]
 \label{lem:netzer-thom}
 Let $\Gamma$ be a group with a finite generating set $S = S^{-1}$.
 Let $c = \sum_g c_g g$ be an element of $\mathbb{R}[\Gamma]$ 
 satisfying $\sum_g c_g = 0$ (i.e. $c$ is in the augmentation ideal) and $c^* = c$.
 Let $D>0$ be an integer such that if $c_g \neq 0$, then 
 $g$ is a product of at most $2^D$ elements from $S$. 
 Then
 \[
 c + 2^{2 D-1} ||c||_1 \cdot \Delta \in 
 \left\{ \sum_{i=1}^n b_i^* b_i \mid b_1, \dots, b_n \in \mathbb{R}[\Gamma] \right\},
 \]
 where $||c||_1 = \sum_g |c_g|$.
 Moreover, 
 if $S$ does not contain self-inverse elements, then even
 \[
 c + 2^{2 D-2} ||c||_1 \cdot \Delta \in 
 \left\{ \sum_{i=1}^n b_i^* b_i \mid b_1, \dots, b_n \in \mathbb{R}[\Gamma] \right\}.
 \]
 \end{lemma}

Using this lemma we explain how we possibly obtain an exact 
solution from a numerical solution.
The following lemma is implicit in \cite{netzer-thom}.

\begin{lemma}
\label{lem:netzer-thom2}
Suppose a constant $\varepsilon>0$ and 
$b_1, \cdots, b_m \in \mathbb{R}[\Gamma]$ are given. 
Assume  $\sum_{i=1}^m b_i^* b_i$ is in the augmentation ideal. 
Set 
$$c = \Delta^2 - \varepsilon \Delta - \sum_{i=1}^m b_i^* b_i.$$
Suppose $c = \sum_g c_g g$.
Let $D>0$ be an integer such that if $c_g \neq 0$, then 
 $g$ is a product of at most $2^D$ elements from $S$. 
Assume $\varepsilon - 2^{2 D-1} ||c||_1 >0$.
Then, $\Gamma$ is a Kazhdan group and $\varepsilon - 2^{2 D-1} ||c||_1 $
is a spectral gap. 
\end{lemma}

\proof
Notice that $c$ is in the augmentation ideal. 
Since $c^*=c$, we can 
apply Lemma \ref{lem:netzer-thom}
to $c$, and 
 there must exist $b_{m+1}, \dots, b_{m+m'}$ such that
 $c+ 2^{2 D-1} ||c||_1 \cdot \Delta =\sum_{i=m+1}^{m'} b_i^* b_i.$
 Plug in the definition of $c$ to this, and get 
 
\[
\Delta^2 - (\varepsilon - 2^{2 D-1} ||c||_1) \, \Delta = \sum_{i=1}^{m+m'} b_i^* b_i.
\] 
This implies the conclusion. 
\qed

Now once we obtain a numerical solution $\varepsilon, b_i$, then 
we apply Lemma \ref{lem:netzer-thom2}to the solution after we modify
$b_i$ (usually, slightly) so that $c$ is in the augmentation ideal.
Notice that we have a more chance 
to have $\varepsilon - 2^{2 D-1} ||c||_1 >0$ if 
$\varepsilon$ is larger and $D$ and $||c||_1$ are small. 
The moreover part of Lemma \ref{lem:netzer-thom} could be 
used to improve the estimate, but it will not be so critical 
in our experiments. 

\subsection{Semidefinite programming}
Although we can apply Lemma 
\ref{lem:netzer-thom2} to any tuples of $\varepsilon, b_1, \cdots, b_m$, 
we have a better chance to succeed if we start with a good 
numerical solution. 
We briefly explain how we find  a solution by a computer.
Following \cite{netzer-thom} we use {\it Semidefinite programming} (SDP)
to find a solution. We refer interested readers to their paper for details. 

Here we explain only the point that is important for us.
To set up and solve an optimization problem by SDP,
we first fix a positive integer $d$ that is an upper bound 
of the {\it support range} of solutions $b_i$'s of the equation
\ref{eq:ozawa}, namely, the word length 
of the group elements in $b_i$'s that have non-zero
coefficients is at most $d$.
For each $d >0$, we solve an optimization problem
to maximize $\varepsilon \ge 0$ such that a solution in the support range
$\le d$ exists for the equation
\ref{eq:ozawa}.

There are three reasons why the algorithm for a given $d$ does not stop:
(i) a group is not a Kazhdan group; (ii) it is a Kazhdan group but 
a solution with the support range $\le d$ does not exist with $\varepsilon >0$; or (iii)
a solution exists in that support range but 
the computer does not have enough power
to find a solution. 
The bigger the support range $d$ is, the more chance  there is that 
the computation will not be finished. 
Also, even if the algorithm stops for some $d$ and gives a positive $\varepsilon$, 
maybe there is a solution for bigger $d$ that gives a larger $\varepsilon$, 
so that our $\varepsilon$ is smaller than the spectral gap. 
See Section \ref{sec:SL_F} for an example of the numerical result
that shows different bounds for different $d$.
Interestingly, it turns out that sometimes, the $\varepsilon$ we obtain for a small $d$
is very close to the actual spectral gap (for example, Table \ref{table:CMSZ2_q=2}).

We have a heuristic estimate on $d$.
If the longest relators  in a given presentation of a group have the word length 
$< 4\ell$, then our experiments tend to work 
if we take  $d = \ell$.
(Notice that the word length of the group elements
that appear in the equation \ref{eq:ozawa} is at most $2d$.
We then identify two such words as group elements
using a relator of length at most $2d+2d=4d$.
So, $d$ should be at least as big as $\ell$, otherwise
there are relators that are not taken into account.)
For example, in Section \ref{subsec:tA2_triangle},
the relators for the group $G_T$ has length 3, and 
$d=1$ works for the experiments.

Compared to \cite{netzer-thom}  there is one extra ingredient  in our approach.
To obtain a lower bound of the Kazhdan constant of
$\mathrm{SL}(3,\Bbb Z)$, only integer matrices 
and rational numbers appear in the algorithm to solve the equation 
\ref{eq:ozawa}. So, a computation is done  without dealing 
with a presentation of the group. 
But we start with a finite presentation 
of a group (except for $\SL(n,\Bbb Z)$ and finite groups)
 then try  to solve the equation \ref{eq:ozawa} .
For that, as a part of the algorithm, we (have to) replace a product of generators with another 
product in the right hand side of the equation properly using a relation of the group. 
But we do not use/need the solution of the word problem of the group. 
Notice that once we find a solution of the equation \ref{eq:ozawa}, even if 
without completely solving the word problem on the way, the $\varepsilon$ we get
is  a lower bound for the group. 
Possibly this part of the algorithm (ie, how much we do the replacement)
  may affect its efficiency. But mathematically speaking, 
the idea is elementary and we skip details. 

\subsection{About tables}
The rest of the paper is  the result of our 
 computer experiments on various examples
 of groups.  
 We make tables
 to present our results. It is mostly about 
lower bounds of the Kazhdan constants, $\kappa$. 

We first
numerically find the maximal $\varepsilon$ for which the equation (\ref{eq:ozawa})
has a solution $\{b_i\}$. If $\varepsilon >0$,  then we plug it in to the left hand side of the inequality (\ref{eq:ek}). In this way we obtain the ``numerical (lower) bound''
of $\kappa$. 

If the support range of $\{b_i \}$ 
is $d$, then we choose $D$
such that $2d \le 2^D$. Then we can apply Lemma \ref{lem:netzer-thom2} to the $\varepsilon>0$ and  the solutions $\{b_i\}$, and obtain 
a ``certified'' lower bound, $\varepsilon - 2^{2 D-1} ||c||_1 $, of the spectral gap. If this is positive,
it will give a ``certified (lower) bound'' of $\kappa$
using the inequality (\ref{eq:ek}).

We occasionally mention known lower/upper bounds of $\kappa$, 
also the exact values of $\varepsilon, \sqrt{2\varepsilon/|S|}, \kappa$
to compare with our bounds. 
Also we sometimes mention the support range $d$ for which we find
the numerical solution $\varepsilon, b_i$, as well as $m$ that is 
the number of the elements $b_i$'s. 

\section{Lattices on $\tA2$-buildings}\label{sec:A2_group}
The first examples of our experiment are 
 {\it uniform lattices}
on $\tA2$-buildings, namely,
 finitely generated groups
that act on $\tA2$-buildings by automorphisms, properly 
and cocompactly.
A general reference is  \cite[S5.7]{BdlHV}.

It is known that those groups are Kazhdan groups.
Moreover,  the Kazhdan constants
are known for a certain family of groups with certain 
generating sets. As far as we know 
this is the only case where the Kazhdan constants
are known for infinite groups. 
We will compare them with our lower bounds. 

\subsection{$\tA2$-buildings}

Let $\Pi$ be a finite projective plane of order $q \ge 2$ with a 
set $P$ of points, a set $L$ of lines
and an incidence relation between lines and points.
Each point (line) is incident with $q+1$ lines (points, resp.).
Also $|P|=|L|= q^2+q+1$.

The {\it incident graph} is a bipartite graph whose vertex
set is $P \cup L$ such that there is an edge 
between $p \in P$ and $\ell \in L$ if $p,\ell$ are incident.

An {\it $\tilde A_2$-building} is a $2$-dimensional, 
simply connected, connected simplicial complex such that
the link of any vertex is the incidence graph of a finite projective plane.


The most familiar example of a projective 
plane is the projective plane $PG(2,F)$ over a filed $F$
(those are called {\it Desarguesian plane}) , namely, 
form a 3-dimensional vector space over $F$, 
letting $P,L$ be the sets of 1- and 2-dimensional
subspaces with incidence being inclusion.
If $|F|=q$, then $PG(2,F)$ is also denoted by $PG(2,q)$. 
(cf. \cite{CMS}).

The easiest example is when $F=F_2$.
In this case, the order of $PG(2,2)$ is $2$, so that
the incidence graph
is a bipartite graph with 14 vertices 
and the degree of each vertex is 3.
The graph is called the {\it Heawood graph}.

We quote a theorem. The first assertion is known 
in various forms (Pansu \cite{pansu}, Zuk \cite{Zuk}, Ballmann-\'{S}wi\k{a}tkowski  \cite{BS}).
For the bound, see  \cite[Theorem 5.7.7]{BdlHV}.
A finite symmetric 
generating set $S$ is explicitly given (see (Pvi) in \cite[Section 5.4]{BdlHV}).

\begin{proposition}\label{A2.bound}
Let $X$ be an $\tA2$-building and suppose a discrete group
$G$ acts on $X$ by automorphisms, property and co-compactly. 
Then $G$ is a Kazhdan group. 

Moreover, if the links of the vertices of $X$
are the incidence graphs of the finite projective planes
of order 
 $q$, 
  then there is a finite generating set $S$ of $G$
 such that  $\sqrt{ \frac{2(\sqrt q -1)^2}{(\sqrt q-1)^2 + \sqrt q}}$
 is a Kazhdan constant.
\end{proposition}

In some cases, $\sqrt{ \frac{2(\sqrt q -1)^2}{(\sqrt q-1)^2 + \sqrt q}}$
 is the Kazhdan constant for $S$ (see Corollary \ref{gap.achieve}).


\subsection{$\tA2$-groups and triangle presentation}
\label{subsec:tA2_triangle}
We review one way to construct uniform lattices
of $\tA2$-buildings following 
\cite{CMS}.
Let $\Pi=(P,L)$ be a finite projective plane and $\lambda: P \to L$
be a bijection. A {\it triangle presentation} compatible with $\lambda$
is a set $T$ of triples $(x,y,z), x,y,z \in P$ such that 

(A) given $x,y \in P$, then $(x,y,z) \in T$ for some $z \in P$ if and only 
if $y,\lambda(x)$ are incident;

(B) $(x,y,z) \in T$ implies $(y,z,x) \in T$;

(C) given $x,y \in P$, then $(x,y,z) \in T$ for at most one $z \in P$.

Given a triangle presentation $T$, we define a group $G_T$,
called an  $\tA2$-{\it group},  as follows:
$$G_T= \langle \{a_x\}_{x \in P} \mid a_x a_y a_z =1\,  {\rm if} \, (x,y,z) \in T \rangle.$$
The set of $a_x$ and their inverses (which are labeled by $L$
such that $a_x^{-1}$ is by $\lambda(x) \in L$) is called {\it the set of natural generators}
(from the triangle presentation).

Then the Cayley graph of $G_T$ is (the 1-skeleton of) a 
(``thick'') $\tA2$-building such that the link
of every vertex is the incidence graph of $(P,L)$. So, 
$G_T$ acts on the building properly and 
cocompactly.
Moreover, the set of natural generators will be $S$
in Proposition \ref{A2.bound}.

In \cite{CMSZ2}, they found all triangle presentations (and $\lambda$)
in the case that $(P,L)$ is the projective planes $PG(2,q)$ for $q=2,3$.

In fact a converse holds, \cite{CMSZ}: 
a group is an $\tA2$-group
if it acts freely and transitively on the vertices
of an $\tA2$-building, and 
if it induces a cyclic permutation of the ``type'' of the 
vertices (there are three types for vertices of a $\tA2$-building).

\cite{CMS}  obtained 
the Kazhdan constant for a class of  $\tA2$-groups
with respect to the natural generators.

\begin{theorem}\label{A2group.bound}\cite[Th 4.6]{CMS}
Let $G$ be an $\tA2$-group obtained from $\Pi=PG(2,q)$.
Let $S$ be the set of natural generators.
Then 
$$\kappa(G,S)= \sqrt{2 \varepsilon_q},$$
where
$$\varepsilon_q= 1- \frac{q(\sqrt q + \sqrt {q^{-1}} +1)}{q^2+q+1}.$$
\end{theorem}

The constant $\sqrt{2 \varepsilon_q}$ is equal to the constant 
in Proposition \ref{A2.bound}. Indeed one can use the proposition 
to show $\sqrt{2 \varepsilon_q} \le \kappa(G,S)$ by checking 
$S$ satisfies the condition of the proposition, but the 
other inequality is hard.

We remark that for  the above $G$ and $S$, we have 
$$\varepsilon(G,S))/|S|=\varepsilon_q,$$
where $\varepsilon(G,S)$ is the spectral gap. 
Indeed, we have $\varepsilon(G,S)/|S| \le \varepsilon_q$
from  Theorem \ref{A2group.bound} and 
Theorem \ref{thm:gap}.
To see the other inequality we first define a graph.
Let $\mathcal G(S)$ be the graph whose vertex set is $S$ such that 
we join $s,t \in S$ if there is $u \in S$ with $s=tu$.
Assume that $\mathcal{G}(S)$ is connected.
Let $\lambda$ be the first positive eigenvalue
of the Laplacian on $\mathcal G(S)$.
Note that in our case, 
 $\mathcal G(S)$
is the link of a vertex of the building, which 
is the incidence graph of $\Pi$. But it is the incidence graph of the finite projective plane of order $q$. So, $\lambda= 1- \sqrt q /(q+1)$, 
see  \cite[Proposition 5.7.6]{BdlHV}.

Now it is pointed out in 
\cite[Example 5]{ozawa} (for the normalized Laplacian), we have 
(for this we only need that $\mathcal{G}(S)$ is connected)
$$(2-\lambda^{-1})|S| \le \varepsilon(G,S).$$
Indeed, he explicitly gives solutions $b_i, i \in S$ of the equation 
\ref{eq:ozawa} for $\varepsilon = (2-\lambda^{-1})|S|$, 
so $(2-\lambda^{-1})|S| \le \varepsilon(G,S)$. 
Moreover, the support of $b_i$ is contained in $S$.
By computation,  $2-\lambda^{-1} = \varepsilon_q$, so
we obtain $\varepsilon_q |S| \le \varepsilon(G,S)$. We are done. 
We have also shown that the equation \ref{eq:ozawa}
has a solution when  $\varepsilon$ is equal to $\varepsilon(G,S)$.


We record this discussion. 
\begin{corollary}\label{gap.achieve}
Let $G$ and $S$ be as in Theorem \ref{A2group.bound}, 
and $\Delta$ the Laplacian.
Then the spectral gap, $\varepsilon(G,S)$, (in the sense of Ozawa)
is achieved, ie, there is a solution $b_i, i \in S$ for the equation \ref{eq:ozawa}
for the $\varepsilon(G,S)$. Moreover, the support of $b_i$ is 
contained in $S$.

The set $S$ satisfies the condition for $S$ in Proposition 
\ref{A2.bound}, and the constant in the proposition is actually
the Kazhdan constant for $S$.

Also, for the $S$ and $\varepsilon(G,S)$,
the equality holds in the inequality \ref{eq:ek}
in Theorem \ref{thm:gap}.
\end{corollary}

Note that Theorem \ref{A2group.bound} is only for the natural generators
from triangle presentations. 
As far as we know, this is the only class
of infinite groups whose Kazhdan constants
are known. 

In view of  Corollary \ref{gap.achieve}, we are curious to know if 
the infimum in the definition of $\kappa(\Gamma,S,\pi), \kappa(\Gamma,S)$
is achieved for the groups in Theorem \ref{A2group.bound}.

\subsection{Computation of $\varepsilon$ of $\tA2$-groups}
As we mentioned, in the case that $q=2,3$, there is a list
of all $\tA2$-groups with the natural generators, \cite{CMSZ2}.

{\it The case $q=2$}.\\
There are 9 presentations: 
$A_1, A_1', A_2, A_3, A_4, B_1, B_2, B_3, C_1$.
From the presentations, we compute lower bounds of the spectral gaps.
See Table \ref{table:CMSZ2_q=2}.
We know $\kappa = 0.465175...$ by letting $q=2$
in Theorem \ref{A2group.bound}.
We observe that our numerical bound is (almost)
identical to the actual value, and 
the certified bound is also very close. 
Also, in view of Corollary \ref{gap.achieve},
we expect to find solution with $d=1$, ie, 
the support of $b_i$ is contained in $\{1\} \cup S$,
which happens in our experiment.

\begin{table}
\begin{tabular}{c|c|c||c}
&  certified bound & numerical solution& $\kappa$  by \cite{CMS} \\
\hline
$A_1$ & 0.465164... & 0.465175... & 0.465175...  \\
$A_1'$ & 0.465166... & 0.465175... & 0.465175... \\
$A_2$ & 0.465167... & 0.465175... & 0.465175...  \\
$A_3$ & 0.465167... & 0.465175... & 0.465175...  \\
$A_4$ & 0.465165... & 0.465175... & 0.465175...  \\
$B_1$ & 0.465164... & 0.465175... & 0.465175...  \\
$B_2$ & 0.465167... & 0.465175... & 0.465175...  \\
$B_3$ & 0.465167... & 0.465175... & 0.465175...  \\
$C_1$ & 0.465167... & 0.465175... & 0.465175...  \\
\end{tabular}
\caption{Lower bounds of $\sqrt{2\varepsilon/|S|}$
(therefore $\kappa$) for the groups on the list in \cite[p. 212]{CMSZ2}
with $q=2$. $|S|=14$. 
We set $d = 1$, then $m=7$. 
}
\label{table:CMSZ2_q=2}
\end{table}

{\it The case $q=3$}.\\
This case contains more groups,
and we compute our bounds only for the first
several ones: groups of 1.1 to 1.8 in their classification. 
See Table \ref{table:CMSZ2_q=3} for the result. 
We know $\kappa = 0.687447... $ by letting $q=3$
in Theorem \ref{A2group.bound}.
Again, our numerical bound is (almost)
identical to the actual value, and 
the certified bound is also very close. 

We see those results as a supporting evidence 
that if we obtain a lower bound
of the spectral gaps, maybe
the bound is not so far from the actual 
value.

\begin{table}
\begin{tabular}{c|c|c||c}
&  certified bound & numerical bound & $\kappa$  by \cite{CMS} \\
\hline
$1.1$ & 0.687430... &  0.687447...  & 0.687447...   \\
$1.1'$ & 0.687430... & 0.687447... & 0.687447... \\
$1.2$ & 0.687431... & 0.687447... & 0.687447... \\
$1.3$ & 0.687430... & 0.687447... & 0.687447... \\
$1.4$ & 0.687429... & 0.687447... & 0.687447... \\
$1.5$ & 0.687431... & 0.687447... & 0.687447... \\
$1.6$ & 0.687433... & 0.687447... & 0.687447... \\
$1.7$ & 0.687432... & 0.687447... & 0.687447... \\
$1.8$ & 0.687432... & 0.687447... & 0.687447... \\
\end{tabular}
\caption{Lower bounds of $\sqrt{2\varepsilon/|S|}$ for some groups 
on the list in \cite[pp. 213--222]{CMSZ2} with $q=3$. $|S|=26$.
 ($d = 1$, $m = 13$. )
 }
\label{table:CMSZ2_q=3}
\end{table}

For both $q=2, 3$, 
we found our $\varepsilon$ for the support range $d=1$.
Namely, in the numerical solutions $b_i$, the coefficients
are non-zero only on the generators. 
We believe that our solutions, $\varepsilon$ and $b_i$'s,
are good numerical approximations of the ones
mentioned in \cite[Example 5]{ozawa}.
For example, observe $2m=|S|$ holds for our solution, 
which also holds for the solutions in \cite[Example 5]{ozawa}

We note that we obtain a positive lower bound for the Kazhdan 
constants only using the group presentations. In this sense
it is a new proof that those groups are Kazhdan groups.


\subsection{Triangles of groups}
\label{sec:ronan}
We now discuss another family of  groups that are lattices on $\tA2$-buildings.
Here is a list:
\begin{equation}
\label{eq:ronan}
\begin{split}
G_1&= \langle a,b,c \mid a^3,b^3,c^3, (ab)^2=ba,(bc)^2=cb, (ca)^2=ac \rangle. \\
G_2&= \langle a,b,c \mid  a^3, b^3, c^3, (ab)^2=ba, (bc)^2=cb, (ac)^2=ca \rangle. \\
G_3&= \langle a,b,c \mid a^3,b^3,c^3,(ab)^2=ba,(ac)^2=ca,(c^{-1}b)^2=bc^{-1} \rangle. \\
G_4&= \langle a,b,c \mid a^3,b^3,c^3,(ab)^2=ba, (ac)^2=ca, (bc^{-1})^2=c^{-1}b \rangle.\\
\end{split}
\end{equation}

Each of the four groups acts on a $\tA2$-building
such that the action on the triangles is regular
and the quotient is one triangle (\cite{Ronan}). 
So, it is a lattice, so that a Kazhdan group. 
Since the quotient has three vertices, 
Theorem \ref{A2group.bound} does not apply to 
this action. 
Also, the generating set, $S$ consisting $a,b,c$ and their inverses, 
does not satisfy the condition for $S$ in Proposition \ref{A2.bound}.
(For example, $\mathcal G(S)$ is not connected, 
with  three components).

Those four groups are interesting and studied from various 
viewpoints. 
Here is a list of facts on the four groups:

\begin{enumerate}

\item
\label{item:ronan}
There is a geometric characterization of the four groups, \cite[Theorem 2.5]{Ronan}:
if $\Delta$ is a ``trivalent triangle geometry'' (ie, a 2-dimensional complex
of triangles such that the link of each vertex is the incidence
graph of $PG(2,2)$), admitting a group
$G$ of automorphisms that is regular on the set of triangles, then $G$ is a quotient group of one of $G_i$.

\item
This is the list of all  fundamental groups
of ``triangles of groups'' (see \cite{St}) such that 
the edge groups are $\Bbb Z/3\Bbb Z$, the vertex
groups are the Frobenius group of order 21, 
$\langle a,b \mid a^3, b^7, a b a^{-1}=b^2 \rangle$,
and 
the face group is trivial.

Each $G_i$ acts on a $\tA2$-building 
properly such that the quotient is one triangle,
 \cite[p118,119]{Tits}.

\item
They are automatic groups (S. M. Gersten and H. Short).
They have a common growth function:
$\frac{1+46z+16z^2}{1-8z+16z^2}$, 
\cite[Example 5.1]{Floyd}.

\item\label{A2T}
They have property (T)
since they are lattices of $\tA2$-buildings. 
\item
$G_1$ and $G_3$ are linear, \cite{KMW2}.
$G_2$ and $G_4$ are not arithmetic, \cite{Tits}.
\item
$G_1$ and $G_2$ are perfect, \cite[Prop 1]{KMW}.
\item
$G_3$ and $G_4$ have normal subgroups
of index 3, which are $\Gamma_1$ and $\Gamma_2$, resp.,  as follows:

$\Gamma_1 = \langle s,t,x \mid s^7=t^7=x^7=1, st=x, s^3t^3=x^3 \rangle$,
cf. \cite[S 5.3]{Es}.

$\Gamma_2 = \langle s,t,x \mid  s^7=t^7=x^7=1, st=x^3, s^3t^3=x \rangle $,
which is  not linear, and 
a subgroup of index 3 in  $G_4$, \cite{BCL}. 
So, $G_4$ is not linear. 

\end{enumerate}

This family is interesting for
us because of (\ref{A2T}).
Here are the numerical results
of our computation.
We obtained positive numbers that 
are very close to each other: 

\begin{table}
\begin{tabular}{c|c|c}
&  certified bound  & numerical bound \\
\hline
$G_1$ & 0.239014... & 0.239146... \\
$G_2$ & 0.238700... & 0.239146... \\
$G_3$ & 0.238405... & 0.239146... \\
$G_4$ & 0.238605... & 0.239146... \\
\end{tabular}
\caption{Lower bounds of $\sqrt{2 \varepsilon / |S|}$,
where $|S|=6$, 
of Ronan's groups (\ref{eq:ronan}).  ($d = 2$.)
We obtain the same numerical bounds for $d=3$.
Theorem \ref{A2group.bound}
does not apply to $S$, so $\kappa$ is unknown.}
\label{table:ronan}

\end{table}

\begin{theorem}
The spectral gap of $G_i$ with respect
to $a,b,c$ and their inverses is at least 0.238.
More generally, if $\Delta$ is a trivalent triangle
geometry with $G$ acting as a regular automorphisms
on the set of triangles of $\Delta$, then 
the spectral gap of $G$ w.r.t. the natural three
generators and their inverses is at least 0.238.

\end{theorem}
\proof
The first assertion is clear from our computation.
The second assertion follows from the fact (\ref{item:ronan}), since 
$G$ is a quotient of one of $G_i$, and the spectral gap,
does not decrease.
\qed

Again, our positive bound is obtained only from the presentations,
so that it is a new proof that those groups are Kazhdan groups. 

Our numerical lower bound 0.239146... is obtained 
from a numerical lower bound of $\varepsilon$, which is 
0.171573. 
We suspect that 0.171573 (and maybe 0.239146...  as well) is a 
good approximation of the spectral gaps
(and  maybe the Kazhdan constants) 
for those four groups. In particular, 
it suggests that those four groups have the 
same spectral gap (and Kazhdan constants) for the natural 
generating sets, $a,b,c$ and their inverses. 

It seems there is no lower bound explicitly given for the Kazhdan constants
of those four groups w.r.t. $S=\{a^\pm, b^\pm, c^\pm\}$.
As we said $S$ does not satisfy the condition in Proposition
\ref{A2.bound}. It is possible to obtain 
a generating set $S'$ using the proposition, for which 
we have $\kappa(G_i,S')=0.465175...$. That would
give a lower bound of $\kappa(G_i,S)$, by writing each element 
of $S'$ as a product of elements of $S$, but it 
will be much smaller than 0.465175..., in particular 
smaller than our lower bound.

We point out that our numerical bound for 
$\varepsilon$  is almost identical to $(\sqrt 2-1)^2
=0.17157287525...$, and that the numerical bound
for $\kappa$ is almost identical to 
$(\sqrt 2 -1) / \sqrt 3=0.23914631173...   $.
We do not know any explanation. 
Since our numerical bound is in fact a ``solution'', 
so we believe 
$$\varepsilon =(\sqrt 2-1)^2$$
holds.

\section{$\mathrm{SL}(3,\mathbb{Z})$, $\mathrm{SL}(4,\mathbb{Z})$ and $\St_3(\Bbb Z)$}
\label{sec:SL3Z}

Next we deal with $\mathrm{SL}(n,\Bbb Z)$.
\cite{netzer-thom} obtained a spectral gap for $\mathrm{SL}(3,\mathbb{Z})$.
Let $E_n$ be the set of all $n \times n$ elementary matrices with $\pm 1$ off the diagonal.
This is a symmetric generating set of $\mathrm{SL}(n,\mathbb{Z})$.

An upper bound of $\kappa(\mathrm{SL}(n,\mathbb{Z}), E_n)$ was given by Zuk 
(see \cite[p. 149]{shalom}), and a lower bound by Shalom \cite{shalom}.
Kassabov improved the lower bound in \cite{kassabov}.

\begin{theorem}[Kassabov \cite{kassabov}, Zuk]
For $\mathrm{SL}(n,\mathbb{Z})$ and the symmetric generating set $E_n$,
the (optimal) Kazhdan constant $\kappa(\mathrm{SL}(n,\mathbb{Z}), E_n)$ satisfies 
\[
\frac{1}{42 \sqrt{n} + 860} \leq \kappa(\mathrm{SL}(n,\mathbb{Z}), E_n) 
\leq \sqrt{\frac{2}{n}}.
\]
\end{theorem}

For $\mathrm{SL}(3,\mathbb{Z})$,  Netzer and Thom \cite{netzer-thom} already
obtained a lower bound that is much better
than the previous ones  
by solving a semidefinite programming based on Theorem \ref{thm:ozawa}.
We improve their bound. For $\mathrm{SL}(4,\mathbb{Z})$
we obtain a new bound that is much better than any other known 
bounds. Our algorithm could not find any bound
for $\mathrm{SL}(5,\mathbb{Z})$.
The results are summarized in Table \ref{table:SL3Z_SL4Z}.

\begin{table}[h]
\begin{tabular}{c|c|c|c||c}
& Kassabov & Netzer-Thom & Certified bound & Upper bound by Zuk \\
\hline
$\mathrm{SL}(3,\mathbb{Z})$ & 0.001072...  & 0.1783... & 0.2155... & 0.8164... \\
$\mathrm{SL}(4,\mathbb{Z})$ & 0.001059... & & 0.3285... & 0.7071... \\
\end{tabular}
\caption{Lower bounds of $\kappa(\mathrm{SL}(n,\mathbb{Z}), E_n)$ for $n = 3,4$. 
($d =2$)}

\label{table:SL3Z_SL4Z}
\end{table}

Here we explain some details. 
As in \cite{netzer-thom}, we try to find solutions $\varepsilon$ and $b_i$'s 
of the equation (\ref{eq:ozawa}) by numerical calculation 
for the support range $d=2$ (with respect to $E_n$).
There are $121$ group elements with the word 
length $\le 2$ in $\mathrm{SL}(3,\mathbb{Z})$ and 
$433$ in $\mathrm{SL}(4,\mathbb{Z})$.
For $d=2$, the certified bound $\varepsilon$ is
$0.278648...$ for $\mathrm{SL}(3,\mathbb{Z})$ 
($0.1905$ in \cite{netzer-thom} . Our numerical bound is $0.2804...$),
and $1.29562...$ for $\mathrm{SL}(4,\mathbb{Z})$
(the numerical bound is $1.313...$).
From the certified bounds for $\varepsilon$,
we obtain our bounds for  the Kazhdan constant.

We discuss the Steinberg groups, $\St_n(\Bbb Z)$, 
which are defined as follows:
\[
\langle x_{i j} \, (i, j \in \{1,2, \dots, n\}, i \neq j) \mid 
[x_{i j}, x_{j k}] = x_{i k} (i \neq k), \, [x_{i j}, x_{k l}] = 1 (i \neq l, j \neq k) \rangle.
\]

It is known that for $n \geq 3$,
$\mathrm{St}_n(\mathbb{Z})$ is an extension of 
$\mathrm{SL}(n,\mathbb{Z})$ by $\mathbb{Z}/2\mathbb{Z}$,
(see \cite[Th 10.1]{milnor}), therefore $\mathrm{St}_n(\mathbb{Z})$ is 
a Kazhdan group for $n \ge 3$.
In fact, $\SL(n,\Bbb Z)$ is obtained from $\St_n(\Bbb Z)$
by adding one relation: $(x_{12} x_{21}^{-1} x_{12})^4=1.$
This element has order 2.
The generators $x_{ij}$
are mapped to the natural generators, $E_n$, of $\SL(n,\Bbb Z)$.
So, $\kappa(\St_n(\Bbb Z), \{x_{ij}\})  \le \kappa(\SL(n,\Bbb Z, E_n) .$

By our computation, we obtain a certified lower bound 
as follows:
\begin{equation}\label{steinberg_kappa}
0.171028... \le \kappa(\mathrm{St}_3(\mathbb{Z}), \{x_{ij}\}).
\end{equation}

We only use the presentation of the group, so this gives a new proof that $\mathrm{St}_3(\mathbb{Z})$ is a Kazhdan group. 
But the computation did not finish for $\mathrm{St}_4(\mathbb{Z})$.


\section{Finite reflection groups}

Now we turn our attention to finite groups. 
All finite groups are Kazhdan groups, (cf. \cite{BdlHV}),
but to find  the Kazhdan constants or the spectral gaps are not easy at all
and they are known only for certain families, for example, 
finite cyclic groups \cite{BH}, 
finite Coxeter groups \cite{kassabov_subspace}, with respect to natural generating sets.
In this section, we compare these results and computer calculations.
We also examine some finite groups whose Kazhdan constants are unknown.

\subsection{Coxeter groups}
\label{sec:coxeter}
It is a classical fact that finite irreducible Coxeter groups 
are classified by Dynkin diagrams. 
For Dynkin diagrams $A_n ,B_n, \cdots$, we denote the corresponding Coxeter groups 
by the same symbols by abuse of notation.
We denote their Coxeter generators as $S_{A_n}, S_{B_n}, \cdots$.

The Kazhdan constants are known. 
\begin{theorem}[Kassabov {\cite[Remark 6.3]{kassabov_subspace}}, 
\label{thm:kassabov.k}
Theorem and Proposition 4 of \cite{BH}, 
Theorems 2.1, 2.2 of \cite{bagno} for $A_n$, $B_n$ and $I_2(n)$]
\[
\begin{split}
\kappa(A_n, S_{A_n}) 
&= \sqrt{\frac{24}{(n+1)^3 - (n+1)} }, \\
\kappa(B_n, S_{B_n}) &= \sqrt{\frac{12}{ n (4 - 3\sqrt{2} + 3(\sqrt{2}-1)n + 2n^2)}},  \\
\kappa(D_n, S_{D_n}) &= \sqrt{\frac{12}{ n (n-1) (2n-1) } } \\
\kappa(F_4, S_{F_4}) &= \sqrt{\frac{14-9\sqrt{2}}{34}},  \quad
\kappa(E_6, S_{E_6}) = \sqrt{\frac{1}{39}}, \quad
\kappa(E_7, S_{E_7}) = \sqrt{\frac{4}{399}}, \\
\kappa(E_8, S_{E_8}) &= \sqrt{\frac{1}{310}}, \quad
\kappa(H_3, S_{H_3}) = \sqrt{\frac{124-48\sqrt{5}}{241}}, \quad
\kappa(H_4, S_{H_4}) = \sqrt{\frac{83-36\sqrt{5}}{409}} \\
\kappa(I_2(n), S_{I_2(n)}) &= 2 \sin \left( \frac{\pi}{2 n} \right) 
\quad (n \geq 3, \, \textrm{dihedral group of order $2 n$}) \\
\end{split}
\]
\end{theorem}
Kassabov also computed the spectral gaps of finite Coxeter groups.
\begin{theorem}[Kassabov {\cite[Remark 6.3]{kassabov_subspace}}] 
\label{thm:kassabov.e}
The spectral gap of the unnormalized Laplacian $\Delta$ is
\[\varepsilon=
4 \left( 1 - \cos \frac{\pi}{h} \right)
\]
where $h$ is the Coxeter number.
\end{theorem}
Here $h = n+1, 2n, 2(n-1), 12, 18, 30, 12, 10, 30, m$ 
for $A_n$, $B_n$, $D_n$, $E_6$, $E_7$, $E_8$, $F_4$, $H_3$, $H_4$, $I_2(m)$ respectively.

 Notice that  the left hand side of the inequality (\ref{eq:ek})
is $\sqrt{\frac{4}{n} \left( 1 - \cos \frac{\pi}{h} \right)}$, which is not equal to
the Kazhdan constant. 

We compute lower
bounds of the spectral gaps, then in our usual manner give lower bounds
of the Kazhdan constants.
See Tables \ref{table.coxeterA}, \ref{table.coxeterB}, \ref{table.coxeterD}, \ref{table.coxeter.Ex},   for the results. 
Our numerical bounds for $\sqrt{2 \varepsilon/|S|}$ are identical 
to the actual value (ie. our numerical bound for $\varepsilon$ 
is identical to the actual constant by Kassabov).

\begin{table}[h]
\begin{tabular}{c|l|l|l|l}
 & certified bound  & numerical bound & $\sqrt{2 \varepsilon/|S|}$ 
&  $\kappa$ \\
\hline
$A_{2}$  & 0.99985... & 1.00000... & 1 & 1.00000... \\
$A_{3}$  & 0.62341... & 0.62491... & 0.62491... & 0.63245... \\
$A_{4}$  & 0.43661... & 0.43701... & 0.43701... & 0.44721... \\
$A_{5}$  & 0.32625... & 0.32738... & 0.32738... & 0.33806... \\
$A_{6}$  & 0.25601... & 0.25694... & 0.25694... & 0.26726... \\
$A_{7}$  & 0.20818... & 0.20856... & 0.20856... & 0.21821... \\
$A_{8}$  & 0.17334... & 0.17364... & 0.17364... & 0.18257... \\
$\vdots$ & & & & \\
\end{tabular}
\caption{Coxeter groups $A_n$. Lower bounds of $\sqrt{2 \varepsilon/|S|}$, 
and the known values by Kassabov. 
 ($d = 2$)}
\label{table.coxeterA}
\end{table}

\begin{table}[h]
\begin{tabular}{c|l|l|l|l}
 & certified & expected & $\sqrt{2 \varepsilon/|S|}$ &  $\kappa$ \\
\hline
$B_{2}$  & 0.76482... & 0.76536... & 0.76536... & 0.76536... \\
$B_{3}$  & 0.42163... & 0.42264... & 0.42264... & 0.43147... \\
$B_{4}$  & 0.27464... & 0.27589... & 0.27589... & 0.28580... \\
$B_{5}$  & 0.19718... & 0.19787... & 0.19787... & 0.20707... \\
$B_{6}$  & 0.14872... & 0.15071... & 0.15071... & 0.15889... \\
$B_{7}$  & 0.11558... & 0.11969... & 0.11969... & 0.12689... \\
$B_{8}$  & 0.09053... & 0.09801... & 0.09801... & 0.10437... \\
$\vdots$ & & & & \\
\end{tabular}
\caption{Coxeter groups $B_n$. ($d = 3$)}
\label{table.coxeterB}
\end{table}

\begin{table}[h]
\begin{tabular}{c|l|l|l|l}
 & certified bound & numerical bound  & $\sqrt{2 \varepsilon/|S|}$ &  $\kappa$ \\
\hline
$D_{4}$  & 0.36556... & 0.36602... & 0.36602... & 0.37796... \\
$D_{5}$  & 0.24553... & 0.24677... & 0.24677... & 0.25819... \\
$D_{6}$  & 0.18044... & 0.18063... & 0.18063... & 0.19069... \\
$D_{7}$  & 0.13805... & 0.13953... & 0.13953... & 0.14824... \\
$D_{8}$  & 0.11146... & 0.11196... & 0.11196... & 0.11952... \\
$\vdots$ & & & & \\
\end{tabular}
\caption{Coxeter groups $D_n$.  ($d = 2$)}
\label{table.coxeterD}
\end{table}

\begin{table}[h]
\begin{tabular}{c|l|l|l|l}
 & certified bound  & numerical bound& $\sqrt{2 \varepsilon/|S|}$ &  $\kappa$ \\
\hline
$E_{6}$  & 0.15032... & 0.15071... & 0.15071... & 0.16012... \\
$E_{7}$  & 0.09203... & 0.09317... & 0.09317... & 0.10012... \\
$E_{8}$  & 0.05164... & 0.05233... & 0.05233... & 0.05679... \\
$F_{4}$  & 0.18334... & 0.18459... & 0.18459... & 0.19342... \\
$H_{3}$  & 0.25520... & 0.25545... & 0.25545... & 0.26299... \\
$H_{4}$  & 0.07316... & 0.07401... & 0.07401... & 0.07820... \\
\end{tabular}
\caption{Coxeter groups of exceptional types. 
$d=2$ for $E_n$ (simply laced) and $d = 3$ for the others.}
\label{table.coxeter.Ex}
\end{table}

\subsection{Complex reflection groups}
\label{sec:cpx_ref}
Next we check (finite) complex reflection groups. 
The irreducible (ie, not a product)  ones are classified
into an infinite families $G(m,p,n)$ and 34 exceptional cases.
We apply our algorithms to some of them. 

Let $\mathfrak{S}_n$ be the symmetric group of $n$ elements. 
For $\sigma \in \mathfrak{S}_n$ and $a_1, \dots, a_n \in \mathbb{C}$, we let 
$[(a_1, \dots, a_n), \sigma]$ be the $n \times n$ matrix whose $(i.j)$-entry 
is $a_i$ if $(i,j) = (i,\sigma(i))$ and $0$ otherwise.
For $m,p,n \in \mathbb{N}$ with  $p | m$, let
\[
G(m,p,n) = \{ \, [(a_1, \dots, a_n), \sigma] \mid \sigma \in \mathfrak{S}_n, \, 
a_i \in \mathbb{C}, \, a_i^m=1, \, \biggl( \prod_{j=1}^n a_j \biggr)^{m/p} = 1 \},
\]
which is a finite subgroup of $\mathrm{U}(n)$.

By definition, $G(m,1,n)$ is isomorphic to the wreath product 
$(\mathbb{Z}/ m\mathbb{Z})^n \wr \mathfrak{S}_n$, and $G(m,p,n)$ is an index $p$ subgroup 
of it.
Let $\zeta_m = \exp \left( \frac{2 \pi \sqrt{-1}}{m} \right)$.

The following set generates $G(m,p,n)$.
\begin{itemize}
\item $G(m,1,n)$:
\[
[(\zeta_m, 1, \dots, 1), \mathrm{id}], \quad [(1, \dots, 1), (i,i+1)] \,\,
(i = 1, \dots, n-1)
\]
\item $G(m,m,n)$:
\[
[(\zeta_m^{-1}, \zeta_m, 1, \dots, 1), (1,2)], \quad [(1, \dots, 1), (i,i+1)] \,\,
(i = 1, \dots, n-1)
\]
\item $G(m,p,n)$ $(1 < p < m, \, p | m)$:
\[
\begin{split}
[(\zeta_{m/p}&, 1, \dots, 1), \mathrm{id}], \quad
[(\zeta_m^{-1}, \zeta_m, 1, \dots, 1), (1,2)], \\
&[(1, \dots, 1), (i,i+1)] \,\,
(i = 1, \dots, n-1)
\end{split}
\]
\end{itemize}

Moreover, a group presentation with respect to this generating system 
has been obtained in \cite[Proposition 3.2]{BMR}.
In this paper, we denote this generating system by $S_{G(m,p,n)}$.
As groups with generators, we have
\[
\begin{split}
(G(1,1,n), S_{G(1,1,n)} \setminus \{ 1 \}) \cong (A_{n-1}, S_{A_{n-1}}), \quad 
(G(2,1,n), S_{G(2,1,n)}) \cong (B_n, S_{B_n}), \\
(G(2,2,n), S_{G(2,2,n)}) \cong (D_n, S_{D_n}), \quad 
(G(m,m,2), S_{G(m,m,2)}) \cong (I_2(m), S_{I_2(m)}).
\end{split}
\]

It seems the Kazhdan constant is unknown for 
complex reflection groups, but 
the Kazhdan constants, $\hat\kappa$,
for $G(m,1,n)$ for irreducible representations
(namely, in the definition we only look at irreducible representations)
are known.  
Define 
\[ 
\hat{\kappa}(\Gamma, S) = \inf \{ \kappa( \Gamma, S, \pi) \mid 
\textrm{$\pi$ is irreducible} \}.
\]

From the definition, we have $\kappa \le \hat\kappa$.

\begin{theorem}[\cite{bagno}]
\[
\hat{\kappa}(G(m,1,n), S_{G(m,1,n)}) = 
\sqrt{\frac{|1-\zeta_m|^2}{\displaystyle\sum_{j=1}^n \left(1 + \frac{|1-\zeta_m|}{\sqrt{2}}(j-1)\right)^2}}
\]
\end{theorem}

Now, here is the result of our computation
for $G(m,1,n)$, Table \ref{table:G_m_1_n}.
For $G(m,1,n)$ there is a gap
between our certified/numerical bound and $\hat\kappa$, 
so that we suspect the equality does not 
hold in inequality \ref{eq:ek}, at least for those families. 

\begin{table}[h]
\begin{tabular}{c|l|l|l}
 & certified bound & numerical bound & $\hat{\kappa}$ by Bagno \\
\hline
$G(3,1,2)$  & 0.68040... & 0.68177... & 0.71010... \\
$G(3,1,3)$  & 0.38644... & 0.38753... & 0.40997... \\
$G(3,1,4)$  & 0.25605... & 0.25749... & 0.27490... \\
$G(3,1,5)$  & 0.18544... & 0.18685... & 0.20067... \\
$G(3,1,6)$  & 0.14080... & 0.14352... & 0.15476... \\
$G(3,1,7)$  & 0.10850... & 0.11469... & 0.12405... \\
$\vdots$ & & & \\
$G(4,1,2)$  & 0.62387... & 0.62491... & 0.63245... \\
$G(4,1,3)$  & 0.36392... & 0.36602... & 0.37796... \\
$G(4,1,4)$  & 0.24449... & 0.24677... & 0.25819... \\
$G(4,1,5)$  & 0.17895... & 0.18063... & 0.19069... \\
$\vdots$ & & & \\
$G(5,1,2)$  & 0.55741... & 0.56301... & 0.56341... \\
$G(5,1,3)$  & 0.34022... & 0.34078... & 0.34752... \\
$G(5,1,4)$  & 0.23243... & 0.23374... & 0.24173... \\
$\vdots$ & & & \\
$G(6,1,2)$  & 0.50135... & 0.50462... & 0.50544... \\
$G(6,1,3)$  & 0.31341... & 0.31469... & 0.32037... \\
$G(6,1,4)$  & 0.21875... & 0.21964... & 0.22654... \\
$\vdots$ & & & \\
\end{tabular}
\caption{Lower bounds for $\sqrt{2 \varepsilon / |S|}$ and $\hat\kappa$ for $G(m,1,n)$.
($d = 3$)}
\label{table:G_m_1_n}
\end{table}

Next, we check $G(m,m,n)$.
To compare with our numerical results, 
we give an upper bound of $\kappa(G(m,m,n), S_{G(m,m,n)})$.
\begin{proposition}[Upper bound of $\kappa$]
\label{upper.bound}
For $m \geq 2$ and $n \geq 2$, we have
\[
\kappa(G(m,m,n), S_{G(m,m,n)}) 
\leq \sqrt{\frac{2|1 - \zeta_{2m}|^2}{2+\sum_{j=1}^{n-2} 
\left| 1+|1-\zeta_{2m}| j \right|^2 }}
\]
\end{proposition}
We remark that when $n=2$, although $(G(m,m,2), S_{G(m,m,2)}) \cong (D_n,S_{D_n})$, 
the  bound does not coincide with $\kappa(D_n,S_{D_n})$.
\begin{proof}
If we let
\[
\eta = ( 1, \,\,  \zeta_{2m}, \,\, \zeta_{2m} + \zeta_{2m} |1 - \zeta_{2m}|, \,\, \dots, \,\,
\zeta_{2m} + \zeta_{2m} |1 - \zeta_{2m}|(n-2)),
\]
then $||\eta||^2 = 2+\sum_{j=1}^{n-2} \left| 1+|1-\zeta_{2m}| j \right|^2$
and $||s \cdot \eta -\eta||^2 = 2|1 - \zeta_{2m}|^2$ for all $s \in S_{G(m,m,n)}$.
So we obtain the desired upper bound.
\end{proof}

Our numerical results for $G(m,m,n)$
are in Table
\ref{table:G_3_3_n}

\begin{table}[h]
\begin{tabular}{c|l|l|l}
 & certified bound & numerical bound& upper bound in 
 Prop \ref{upper.bound}  \\
\hline
$G(3,3,2)$  & 0.99999... & 1.00000... & 1.00000... \\
$G(3,3,3)$  & 0.55621... & 0.55851... & 0.70710... \\
$G(3,3,4)$  & 0.34256... & 0.34350... & 0.53452... \\
$G(3,3,5)$  & 0.23364... & 0.23596... & 0.42640... \\
$G(3,3,6)$  & 0.16757... & 0.17445... & 0.35355... \\
$\vdots$ & & & \\

\hline

$G(4,4,2)$  & 0.76482... & 0.76536... & 0.76536... \\
$G(4,4,3)$  & 0.46141... & 0.46538... & 0.55780... \\
$G(4,4,4)$  & 0.30330... & 0.30528... & 0.43136... \\
$\vdots$ & & & \\

\hline
$G(5,5,2)$  & 0.61594... & 0.61803... & 0.61803... \\
$G(5,5,3)$  & 0.38720... & 0.38968... & 0.45950... \\
$G(5,5,4)$  & 0.26542... & 0.26714... & 0.36124... \\
$\vdots$ & & & \\

\end{tabular}
\caption{Lower bounds of $\sqrt{2 \varepsilon / |S|}$ for $G(3,3,n)$;
$G(4,4,n)$.;
and $G(5,5,n)$.
($d = 3$ for $n = 2$, $d = 4$ for $n = 3, 4$ ).
}
\label{table:G_3_3_n}
\end{table}



We check a few more families for $G(m,p,n)$, which 
is in Table \ref{table:G_4_2_n}.


\begin{table}[h]
\begin{tabular}{c|l|l}
 & certified bound & numerical bound \\
\hline
$G(4,2,2)$  & 0.91909... & 0.91940... \\
$G(4,2,3)$  & 0.49128... & 0.49288... \\
$G(4,2,4)$  & 0.30769... & 0.30883... \\
$G(4,2,5)$  & 0.21409... & 0.21585... \\
$\vdots$ & & \\

\hline
$G(6,2,2)$  & 0.77478... & 0.77740... \\
$G(6,2,3)$  & 0.42686... & 0.43170... \\
$G(6,2,4)$  & 0.27576... & 0.27775... \\
$\vdots$ & & \\

\hline
$G(6,3,2)$  & 0.81608... & 0.81649... \\
$G(6,3,3)$  & 0.45887... & 0.46055... \\
$G(6,3,4)$  & 0.29070... & 0.29538... \\
$\vdots$ & & \\

\end{tabular}
\caption{Lower bounds of $\sqrt{2 \varepsilon / |S|}$ for $G(4,2,n)$.
($d = 3$);
$G(6,2,n)$ ($d = 3$);
$G(6,3,n)$ ($d = 3$ for $n = 2$, $d = 4$ for $n = 3, 4$).
}
\label{table:G_4_2_n}
\end{table}



\section{$\mathrm{SL}(n,\mathbb{F}_p)$}
\label{sec:SL_F} 
Our last examples are $\mathrm{SL}(n,\mathbb{F}_p)$.
The Kazhdan constants are not known. 
But the following is known.

\begin{theorem}[Kassabov {\cite[Theorem A'']{kassabov}}]
\label{thm:kassabov_finite_SL}
The Kazhdan constant for $\mathrm{SL}(n,\mathbb{F}_p)$ with respect to the set 
$E_n$ of elementary matrices with $\pm 1$ off the diagonal satisfies
\begin{equation}
\label{eq:kassabov_finite_SL}
\kappa (\mathrm{SL}(n, \mathbb{F}_p), E_n ) \geq  \frac{1}{31 \sqrt{n} + 700}.
\end{equation}
\end{theorem}

We compute lower bounds
for $\mathrm{SL}(n,\mathbb{F}_p)$ with respect to the set $E_n$ 
with $n,p$ small.
Using this example, we give an idea on how
our bounds possibly depend on the support range $d$.  Table \ref{SL_d=2}
is for the support range $d=2$, and 
Table \ref{SL_d=3} is for the support range $d=3$.
We may obtain a better (ie, larger) bound
for a larger $d$, but there is more chance that the computation 
does not finish. 
For example, the bound for $\mathrm{SL}(2,\mathbb{F}_5)$
improves much if we change $d=2$ to $d=3$.

\begin{table}[h!]
\begin{tabular}{c|lll|l}
$p$ & 3 & 5 & 7 & $\kappa_0$ \\
\hline
$\mathrm{SL}(2,\mathbb{F}_p)$ & 
0.7961... & 0.2580... &  ($\ast 1$) & 0.001344... \\
$\mathrm{SL}(3,\mathbb{F}_p)$ &
0.6716... & 0.4981... & 0.3508... & 0.001326... \\ 
$\mathrm{SL}(4,\mathbb{F}_p)$ &
0.5974... & 0.4812... & 0.3284... & 0.001312... \\ 
\end{tabular}
\caption{Certified lower bounds of the Kazhdan 
constant ($d= 2$), and Kassabov's lower bound $\kappa_0$, 
the right hand side of the inequality (\ref{eq:kassabov_finite_SL}).
 Our numerical bound for $\mathrm{SL}(2,\mathbb{F}_7)$ is close to $0$, so that we could not get a positive 
number as a certified bound $^{(\ast 1)}$.
}
\label{SL_d=2}
\end{table}

\begin{table}[h]
\begin{tabular}{c|lll|l}
$p$ & 3 & 5 & 7  & $\kappa_0$ \\
\hline
$\mathrm{SL}(2,\mathbb{F}_p)$ &
0.7958... & 0.6145... & 0.5387... & 0.001344... \\
$\mathrm{SL}(3,\mathbb{F}_p)$ & 
0.6683... & 0.4410... & ($\ast 2$) & 0.001326... \\ 
\end{tabular}
\caption{Certified lower bounds of the Kazhdan 
constant. ($d= 3$.)  The computation did not finish for $\mathrm{SL}(3,\mathbb{F}_7)$ $^{(\ast 2)}$.}
\label{SL_d=3}
\end{table}

\section{Problems suggested by the experimental results}
We mention problems that naturally arise from our experiments.

\begin{enumerate}
\item
Is our bound, 0.239, for the Kazhdan constants
of $G_1, G_2,G_3,G_4$ in Section 3.3 close to the actual values?
Do the four groups have same Kazhdan constants?
Is the Kazhdan constant equal to 
$(\sqrt 2 -1) / \sqrt 3)=0.23914631173...   $ ?

Similarly, is our bound 0.171573 close to 
the spectral gaps $\varepsilon$ of the four groups?
Do they have the same spectral gaps?
Is it $(\sqrt 2 -1 )^2=0.17157287525...$?

To answer the first question for $G_i$, we only need to find 
one representation $\pi$ and one
vector $\xi \in \mathcal H$ such that $\max_{s\in S} ||\pi(s) \xi - \xi||$
is close to  0.239, where $S$ consists of $a,b,c$ and the inverses. 
If one finds such representation, it also implies that 
0.171573 is close to $\varepsilon$.

The numerical solutions $b_i$ we find for the above spectral 
gap are for $d=2$, ie,
the support of $b_i$ are contained in $\{1\} \cup S \cup S^2$.
Is there any explanation for that? (cf. Corollary \ref{gap.achieve},
where we find a solution for $d=1$.)

\item
Among $\mathrm{SL}(n,\Bbb Z), n \ge 3$, 
is the spectral gap $\varepsilon$ 
 a monotone decreasing function on $n$ 
 w.r.t. the generating sets of elementary matrices ?
 How about the Kazhdan constants ?
 Our lower bound for $\mathrm{SL}(4,\Bbb Z)$
 is larger than the one for $\mathrm{SL}(3,\Bbb Z)$.
If it is not monotone, 
for which $n$, are the Kazhdan constant of $\mathrm{SL}(n,\Bbb Z)$
 largest? Combining our lower bound for $\mathrm{SL}(4,\Bbb Z)$
 and Zuk's upper bound,
 it must be at most $18$.

\item
Does the equality hold in the inequality \ref{eq:ek}
for finite complex reflection groups ?
In view of this problem,
compute the spectral gaps (maybe using the method by Kassabov, 
see \cite[Remark 6.5]{kassabov_subspace}).

We believe that our (in particular, numerical) bounds for the spectral gap
are very close to the actual values, and that
our results suggest that the equality does not hold.

\item
What is the behavior of the spectral 
gap (or the Kazhdan constant)
of $\mathrm{SL}(n,\mathbb{F}_q)$ when $q$ increases with $n$ fixed ?
Does it converge to the spectral gap (or the Kazhdan constant) of $\mathrm{SL}(n, \Bbb Z)$
as $q \to \infty$?
Notice that the answer is negative for the Kazhdan 
constant for $n=2$. This is because 
$\mathrm{SL}(2,\Bbb Z)$ does not have property (T),
so that $\kappa=0$, while there is a uniform positive
lower bound for $\mathrm{SL}(2, \Bbb F_p, E_2)$ by Kassabov
(see the discussion around the end of the introduction 
in \cite{BdlHV}).

How about  when $n$ increases with  $q$ fixed ?
Our experiments may suggest that they are monotone on $q$
with $n$ fixed;  and also on $n$ with $q$ fixed.

\item
For a given $\varepsilon < \varepsilon(G,S)$,
the equation \ref{eq:ozawa} has solutions 
$b_i$ (\cite{ozawa}). But is there a solution for $\varepsilon(G,S)$?
See Corollary \ref{gap.achieve}.

Is it possible to estimate the support range $d$
of the solutions in advance? Is there a universal upper bound 
of $d$ for all $\varepsilon$?

\end{enumerate}


\begin{thebibliography}{99}

\bibitem[BH]{BH}
Roland Bacher and Pierre de la Harpe, 
\textit{Exact values of Kazhdan constants for some finite groups}, 
J. Algebra 163 (1994) 495--515.

\bibitem[BCL]{BCL}
Uri Bader, Pierre-Emmanuel Caprace, Jean Lecureux,
\textit{On the linearity of lattices in affine buildings and ergodicity of the singular Cartan flow}, 
arXiv:1608.06265

\bibitem[Bag]{bagno}
Eli Bagno, 
\textit{Kazhdan constants of some colored permutation groups},
Journal of Algebra 282 (2004) 205--231.

\bibitem[B\'{S}]{BS}
W. Ballmann,  J. \'{S}wi\k{a}tkowski, 
\textit{On $L^2$-cohomology and property (T) for automorphism groups of polyhedral cell complexes},
Geom. Funct. Anal. 7 (1997), no. 4, 615--645. 

\bibitem[BdlHV]{BdlHV}
Bachir Bekka, Pierre de la Harpe, and Alain Valette,
\textit{Kazhdan's property (T)}, 
New Mathematical Monographs, 11. Cambridge University Press, Cambridge, 2008. xiv+472 pp.

\bibitem[BMR]{BMR}
Michel Brou\'{e}, Gunter Malle, and Rapha\"{e}l Rouquier,
\textit{Complex reflection groups, braid groups, Hecke algebras},
J. Reine Angew. Math. 500 (1998), 127--190. 

\bibitem[BZ]{brown-ozawa}
Nathanial P. Brown and Narutaka Ozawa, 
\textit{$C^*$-algebras and finite-dimensional approximations}, 
Graduate Studies in Mathematics, vol. 88, 
American Mathematical Society, Providence, RI, 2008.


\bibitem[CMSZ]{CMSZ}
D I Cartwright, A M Mantero, T Steger, A Zappa, 
\textit{Groups acting simply transitively on the vertices of a building of type $\tA2$, I}, 
Geom. Dedicata 47 (1993), no. 2, 143–-166. 

\bibitem[CMSZ2]{CMSZ2}
D I Cartwright, A M Mantero, T Steger, A Zappa, 
\textit{Groups acting simply transitively on the vertices of a building of type $\tA2$, 
II: The cases $q=2$ and $q=3$}, 
Geom. Dedicata 47 (1993), no. 2, 167–-223. 
 
\bibitem[CMS]{CMS}
D I Cartwright, W M{\l}otkowski, T Steger, 
\textit{Property (T) and $\tA2$ groups}, 
Ann. Inst. Fourier (Grenoble) 44 (1994) 213–-248 

\bibitem[Es]{Es}
Jan Essert, 
\textit{A geometric construction of panel-regular lattices for buildings of types $\tA2$ and $\tilde C_2$}, 
Algebr. Geom. Topol. 13 (2013), no. 3, 1531–-1578. 

\bibitem[FP]{Floyd}
William Floyd, Walter Parry,
\textit{The growth of nonpositively curved triangles of groups},
Invent. Math. 129 (1997), no. 2, 289–-359. 

\bibitem[Kas1]{kassabov}
Martin Kassabov, 
\textit{Kazhdan constants for $SL_n(Z)$}, 
Internat. J. Algebra Comput. 15 (2005), no. 5--6, 971--995.

\bibitem[Kas2]{kassabov_subspace} 
M. Kassabov,
\textit{Subspace arrangements and property T}, 
Groups Geom. Dyn. 5 (2011), no. 2, 445--477.

\bibitem[Ka]{Ka}
Ryota Kawakami.
Kazhdan's Property(T) and
Semi-Definite Programming.
Master Thesis paper, 2015, Kyoto University.

\bibitem[KMW]{KMW}
Peter K{\"o}hler, Thomas Meixner, Michael Wester, 
\textit{Triangle groups}, 
Comm. Algebra 12 (1984), no. 13--14, 1595–-1625.


\bibitem[KMW2]{KMW2}
Peter K{\"o}hler, Thomas Meixner, Michael Wester,
\textit{The affine building of type $\tA2$ over a local field of characteristic two},
Arch. Math. (Basel) 42 (1984), no. 5, 400–-407. 

\bibitem[Oz]{ozawa}
Narutaka Ozawa, 
\textit{Noncommutative real algebraic geometry of Kazhdan's property (T)}, 
\href{http://arxiv.org/abs/1312.5431}{arXiv:1312.5431}

\bibitem[Mil]{milnor}
John Milnor, 
\textit{Introduction to algebraic $K$-theory}. 
Annals of Mathematics Studies, No. 72. Princeton University Press, xiii+184 pp. 

\bibitem[NT]{netzer-thom}
Tim Netzer and Andreas Thom,
\textit{Kazhdan's Property (T) via Semidefinite Optimization}, 
\href{http://arxiv.org/abs/1411.2488}{arXiv:1411.2488}

\bibitem[Pan]{pansu}
Pierre Pansu,
\textit{Formules de Matsushima, de Garland et propri\'{e}t\'{e} (T) pour des groupes agissant sur des espaces symétriques ou des immeubles},
Bull. Soc. Math. France 126 (1998), no. 1, 107--139. 

\bibitem[R]{Ronan}
M. A. Ronan,
\textit{Triangle geometries},
J. Combin. Theory Ser. A 37 (1984), no. 3, 294–-319. 

\bibitem[Sha]{shalom}
Yehuda Shalom, 
\textit{Bounded generation and Kazhdan's property $(T)$}, 
Publications Math\'{e}matiques de l'IH\'{E}S 90 (1999), 145--168.

\bibitem[ST]{shephard-todd}
G. Shephard and J. Todd,
\textit{Finite unitary reflection groups},
Canadian J. Math. 6, (1954). 274--304.

\bibitem[Sp]{David Speyer}
David Speyer,
\textit{How feasible is it to prove Kazhdan's property (T) by a computer?},
an answer to MathOverflow question,
\url{http://mathoverflow.net/a/154459} 

\bibitem[St]{St}
John R. Stallings, 
\textit{Non-positively curved triangles of groups}, 
Group theory from a geometrical viewpoint (Trieste, 1990), 491–-503, World Sci. Publ., River Edge, NJ, 1991. 

\bibitem[Ti]{Tits}
Jacques Tits, 
\textit{Buildings and group amalgamations}, 
Proceedings of groups—St. Andrews 1985, 110--127,
London Math. Soc. Lecture Note Ser., 121, Cambridge Univ. Press, Cambridge, 1986. 

\bibitem[Zuk]{Zuk}
A. Zuk, 
\textit{La propri\'{e}t\'{e} (T) de Kazhdan pour les groupes agissant sur les poly\`{e}dres}, 
C. R. Acad. Sci. Paris 323, Serie I (1996), 453--458.




\end{thebibliography}
\end{document}